\documentclass{article}
\usepackage{threeparttable}
\usepackage{index}
\usepackage{authblk}
\usepackage{subfigure}
\usepackage{multirow}
\usepackage{graphicx}
\usepackage{algorithm}
\usepackage{algorithmic}
\usepackage[a4paper,hmargin={2.54cm,2.54cm},vmargin={3.17cm,3.17cm}]{geometry}
\usepackage{amsmath,amssymb,amsthm,cases}
\usepackage{cite}
\usepackage{enumerate}
\usepackage[pdfborder={0 0 0},colorlinks=true,linkcolor=blue,CJKbookmarks=true]{hyperref}

\usepackage{xcolor}
\usepackage{tcolorbox}
\usepackage[normalem]{ulem}
\usepackage{soul}
\usepackage{amssymb}
\usepackage{pifont}

\newtheoremstyle{mystyle}{3pt}{3pt}{}{0cm}{}{}{1em}{}
\theoremstyle{mystyle}
\newtheorem{definition}{\textbf{Definition}}
\newtheorem{theorem}{\textbf{Theorem}}

\newtheorem{lemma}{\textbf{Lemma}}
\newtheorem{proposition}{ \textbf{Proposition}}
\newtheorem{assumption}{ \textbf{Assumption}}
\newtheorem{corollary}{ \textbf{Corollary}}
\newtheorem{remark}{\textbf{Remark}}

\newcommand{\keywords}[1]{\textbf{Keywords:}\quad #1}

\allowdisplaybreaks[4]

\title{\textbf{Inertial Bregman Proximal Gradient Algorithm For Nonconvex Problem with Smooth Adaptable Property}}

\author[1]{Xiaoya Zhang,
  Hui Zhang\thanks{Email: \texttt{h.zhang1984@163.com} },
  Wei Peng
}

\affil[1]{ Department of Mathematics, National University of Defense Technology,
Changsha, 410073, Hunan, China.}

\begin{document}

\maketitle

\begin{abstract}
In this paper we study the problems of minimizing the sum of two nonconvex functions: one is differentiable and satisfies smooth adaptable property. The smooth adaptable property, also named relatively smooth condition, is weaker than the globally gradient Lipschitz continuity. 
We analyze an inertial version of the Bregman Proximal Gradient (BPG) algorithm and prove its stationary convergence. Besides, we prove a sublinear convergence of the inertial algorithm. Moreover, if the objective function satisfies Kurdyka--{\L}ojasiewicz (KL) property,  its global convergence to a critical point of the objective function can be also guaranteed. 
\end{abstract}



\section{Introduction}
In this paper, we consider minimizing the sum of two functions: a proper continuously differentiable function $f$(not necessarily convex) and a proper lower-continuous function(not necessarily convex) $g$: 
\begin{align*} \label{Eq:P}
\inf \{ \Psi(x) := f(x) + g(x): x \in \mathbb{R}^d \} \tag{P}.
\end{align*}
Problem (P) arises in many applications including compressed sensing\cite{donoho2006compressed}, DC-programming for digital communication system\cite{alvarado2014new}, signal recovery\cite{beck2013sparsity}, phase retrieve problem \cite{luke2017phase}.

Although First-Order-Methods for solving convex problems have long history, relative  
algorithms and analysis for totally nonconvex and nonsmooth problems, such as (P), are new\cite{bolte2014proximal}\cite{allen2016variance}\cite{reddi2016proximal}. One common algorithm for solving (P) is based on computing a proximal operator during each step (under certain assumptions which is used to guarantee the existence of the minimizer $x^{k+1}$):
$$ x^{k+1} = \arg\min_{x \in \mathbb{R}^d} \{ g(x) + \langle \nabla f(x^k), x-x^k \rangle + \frac{1}{2\lambda_k}\|x-x^k\|^2 \}.$$
But these nonconvex algorithms have restriction: the gradient of the smooth part $f$ has to be globally Lipschitz continuous on $\mathbb{R}^d$, like \cite{lu2014generalized}\cite{reddi2016proximal}\cite{peng2018nonconvex}, etc. 

Recently,  as to convex problem, Bauschke, Bolte and Teboulle (BBT)\cite{bauschke2016descent} solved this longstanding issue and avoiding globally Lipschitz continuous gradient by leading into a new definition-$L$-smooth adaptable($L$-smad) property. A general Bregman Proximal Gradient(BPG) algorithm who used this property to instead gradient Lipschitz continuity is followed:
$$ x^{k+1} = \arg\min_{x \in \mathbb{R}^d} \{g(x) + \langle \nabla f(x^k), x-x^k \rangle + \frac{1}{\lambda_k}D_h(x, x^k) \}.$$
with Bregman distances instead of the quadratic terms.
Later, in \cite{bolte2017first}, the authors analyzed the nonconvex case and proposed nonconvex BPG algorithm.
And in \cite{teboulle2018simplified}, Teboulle made a review and pointed out that accelerations for this Bregman Proximal Gradient algorithms also have faster rate according to numerical experiments, but it lacks theoretical support.  

In this paper, we proposed an inertial version for Bregman Proximal Gradient algorithm. Inertial algorithms are a focus of optimization, which has been well explored by related researchers. Inertial method was explained as an explicit finite differences discretization of the so-called Heavy-ball with friction dynamical system
(where $g(x)\equiv 0$), see from \cite{ochs2014ipiano}. However, in \cite{ochs2014ipiano} it required that $g$ is convex, we omit this requirement in our paper. This means our algorithm is more general than \cite{ochs2014ipiano}.
Under the same conditions in \cite{bolte2017first}, similar convergence results in \cite{li2015accelerated} can be get for this general framework. It can deal with many nonconvex nonsmooth problems who has no globally  gradient Lipschitz continuity, like quadratic inverse problem, which was not solved by other inertial algorithms, such as in \cite{ochs2018local}\cite{ochs2014ipiano}.

\textbf{Outline and Contributions.} The paper is organized as follows. We first introduce 
some basic definitions in the first part of Section 2, and Bregman distance and $L$-smooth adaptive property later in Section 2; 
Description for inertial Bregman Gradient Descent(iBPG) algorithm follows in the third part. 
A general descent lemma is displayed at the beginning of Section 4.
Convergence analysis of iBPG algorithm for this non-Lipschitz-continuous nonconvex  problem (P) under natural assumptions is analyzed in Section 4.1, which guarantee that any cluster(limiting) point is a critical point.  Besides, a sublinear rate is shown.  Secondly, with additional KL property, we prove that the whole sequence generated by iBPG converges to a critical point in Section 4.2.

\section{Preliminaries}
Throughout the paper, let $\mathbb{N}:= \{0, 1, 2, \dots\}$ be the set of nonnegative integers. We will always work in a finite dimensional Euclidean vector space $\mathbb{R}^d$. The standard Euclidean inner product and and the induced norm and the induced norm on on $\mathbb{R}^d$ are denoted by $\langle \cdot,\cdot\rangle$ and $\| \cdot\|$, respectively. 
Recall some basic notions, $B_{\rho}(\tilde{x}):=\{x \in \mathbb{R}^d : \|x-\tilde{x}\|\leq \rho \}$ as the ball of radius $\rho>0$ around $\tilde{x} \in \mathbb{R}^d$; $\text{dist}(x, \mathcal{S}):= \inf_{y \in \mathcal{S}}\|x-y\|$ as the distance from a point $x \in \mathbb{R}^d$ to a nonempty set $\mathcal{S} \subset \mathbb{R}^d$.

The domain of the function $f : \mathbb{R}^d \rightarrow (-\infty,+\infty]$ is defined by $\text{dom}~f = \{x \in \mathbb{R}^d: f(x)< +\infty\}$. We say that $f$ is proper if $\text{dom}~f \neq \emptyset$. Other generalized notions we employ are refered to \cite{rockafellar2015convex}.

In the following, we recall some important basic definitions used in this manuscript. The first one is subdifferential.
\begin{definition}(subdifferentials\cite[Definition 8.3]{rockafellar2009variational}).
The Fr\'{e}chet subdifferential of $f$ at $\bar{x}\in \text{dom}~f$ is the set 
$\hat{\partial} f(\bar{x})$ of elements $v\in \mathbb{R}^d$ such that
$$ \underset{\substack{x \rightarrow \bar{x} \\ x\neq \bar{x}}}{\lim \inf}\frac{f(x)-f(\bar{x}) - \langle v,x-\bar{x}\rangle}{\|x-\bar{x}\|} \geq 0.$$
For $\bar{x} \notin \text{dom}~f$, we set $\hat{\partial} f(\bar{x}) = \emptyset$. The so-called (limiting) subdifferential of $f$ at $\bar{x} \notin \text{dom}~f$ is defined by
$$ \partial f(\bar{x}):=\{v\in \mathbb{R}^d: \exists x^n \xrightarrow{f} \bar{x}, v^n \in \hat{\partial} f(x^n), v^n \rightarrow v\},$$  
where $x^n \xrightarrow{f} \bar{x}$ means $(x^n,f(x^n)) \rightarrow (\bar{x},f(\bar{x}))$ as $k \rightarrow \infty$, and $\partial f(\bar{x}) = \emptyset$ for $\bar{x} \notin \text{dom}~f$.
\end{definition}

A point $\bar{x} \in \text{dom}~f$ for which $0 \in \partial f(\bar{x})$ is a called a critical point of $f$. 
In this paper, we denote $\text{crit}~\Psi $ as the set of critical points of $\Psi$, which means
$$ \text{crit}~\Psi = \{x: 0 \in \nabla f(x)+\partial g(x) \}.$$

As a direct consequence of the definition of the limiting subdifferential, we have the following closedness property at any $\bar{x} \notin \text{dom}~f$:
$$x^k \xrightarrow{f} \bar{x}, v^k \rightarrow \bar{x}, ~\text{and for all} ~ k \in \mathbb{N}:v^k \in \partial f(x^k) \Longrightarrow \bar{v} \in \partial f(\bar{x}). $$


Finally, we give definition of Kurdyka--{\L}ojasiewicz(KL) property, which was proposed in \cite{bolte2014proximal}. This property would help us to prove some local convergence results. Many functions has KL properties, like semi-algebraic functions，functions definable in an o-minimal structure, and many others\cite{attouch2013convergence}.

\begin{definition}(Kurdyka--{\L}ojasiewicz property\cite{bolte2014proximal})
Let $f : \mathbb{R}^d \rightarrow (-\infty,+\infty]$ be an extended real valued function and let $\bar{x} \in \text{dom}~\partial f$. If there exists $\eta \in [0, \infty]$, a neighborhood $U$ of $\bar{x}$ and a continuous concave function $\psi:[0, \eta) \rightarrow \mathbb{R}_{+}$ such that
$$\psi(0) = 0, \psi \in C^1(0,\eta),~\text{and}~ \psi^{\prime}(s) > 0 ~\text{for~all}~ s \in (0,\eta),$$
and for all $ x \in U \cap [f(\bar{x})<f(x)<f(\bar{x})+\eta]$ the Kurdyka--{\L}ojasiewicz inequality
$$ \psi^{\prime}(f(x)-f(\bar{x}))\text{dist}(0,\partial f(x))\geq 1 $$
holds, then the function has the Kurdyka--{\L}ojasiewicz property at $\bar{x}$.
If, additionally, the function is lsc and the property holds for each point in $\text{dom}~\partial f$, then f is called Kurdyka--{\L}ojasiewicz function.
\end{definition}

\subsection{Smooth Adaptable Functions}
Next, we define the notion of $L$-smooth adaptable property for non-convex $f$ in \cite{bolte2017first}. This property shares the same meaning with the relatively smooth property introduced in \cite{lu2018relatively}. It was extended from the recent work \cite{bauschke2016descent} in which $f$ is convex. This property can analyze functions whose gradient has no global Lipschitz continuous property and function has no convexity itself. It provides a relative smooth notation to a convex function $h$, which is defined as Bregman distance.
We first introduce the definition of Bregman distance\cite{bregman1967relaxation}\cite{bolte2017first}.
\begin{definition}\cite{bolte2017first} 
 Let $S$ be a nonempty, convex and open subset of $\mathbb{R}^d$. Associated with $S$, a function $h:\mathbb{R}^d  \rightarrow (-\infty, \infty] $  is called a kernel generating distance if it satisfies the following:
\begin{enumerate}[(i)]
\item $h$ is proper, lower semicontinuous and convex, with $\text{dom}~h \subset \overline {S}$ and $\text{dom}~\partial h = S$.
\item $h$ is $C^1$ on $\text{int}~\text{dom}~h \equiv S$.
\end{enumerate}
The Bregman distance associated to $h$, denoted by $D_h:\text{dom}~h \times \text{int}~\text{dom}~h \rightarrow [0, +\infty)$ is defined by 
$$D_h(x, y) := h(x) - h(y)- \langle \nabla h(y),x-y \rangle. $$
\end{definition}
Naturally we have that $h$ is convex if and only if 
$D_h(x, y) \geq 0, \forall (x, y) \in \text{dom}~h \times \text{int}~\text{dom}~h.$ 
If in addition $h$ is strictly convex, $D_h(x, y) = 0$ if and only if $x = y$ holds.
Many useful choices for $h$ to generate relevant associated Bregman distances $D_h$ has been listed in many papers. Interested readers could refer to \cite{teboulle2018simplified}.

Throughout the paper we will focus on a smooth $f$ who satisfies the $L$-smooth adaptable condition. It is more general than the condition who describes the gradient Lipschitz continuous property and easier to be satisfied consequently.
Next we will give a clear definition.
\begin{definition}\cite[Definition 2.2.]{bolte2017first} 
A pair $(f, h)$ is called $L$-smooth adaptable($L$-smad) on $C$ if there exists $L > 0$ such that $Lh-g(Lh+g)$ is convex on $S$.
\end{definition}
According to \cite[Lemma 2.1]{bolte2017first}, the pair of functions $(f, h)$ is $L$-smooth adaptable on $S$ if and only if $\|f(x)-f(y)- \langle \nabla f(y), x-y \rangle \| \leq LD_h(x, y)$.
When $h(x)= \frac{1}{2}\|x\|^2$ and consequently $D_h(x,y)=\frac{1}{2}\|x-y\|^2$, the $L$-smad of $f$ would be reduced to gradient Lipschitz continuity.

\section{Inertial Bregman Proximal Gradient Algorithm}
Throughout this paper, we tackle problem (P) under the following assumptions:
\begin{assumption}\label{Assumption:A}
\begin{enumerate}[(i)]
\item $h$ is a kernel generating distance defined in Definition 1 associated with  $S \equiv \mathbb{R}^d$.
\item $f: \mathbb{R}^d \rightarrow (-\infty, \infty]$ is a $C^1$ function and $(f,h)$ is $L$-smad with $\text{dom}~h \subset \text{dom}~f$ on $S \equiv \mathbb{R}^d$.
\item $g: \mathbb{R}^d \rightarrow (-\infty, \infty]$ is  proper, lower semicontinuous.
\item $\Psi^{\ast} := \inf \{ \Psi(u): u \in  \mathbb{R}^d \}􏰂􏱷 > -\infty$.
\end{enumerate}
\end{assumption}

To present an inertial Bregman Proximal Gradient algorithm, we first verify each iteration in is well-posed. It requires another additional assumption for $\Psi$.
\begin{assumption}\label{Assumption:B}
The function $\Psi$ is supercoercive, that is,
$$ \lim_{\|u\|\rightarrow \infty} \frac{\Psi(u)}{\|u\|} = \infty.$$  
\end{assumption}
For all $y \in \text{int}~\text{dom}~h$ and step-size $0< \lambda \leq 1/L$, $0\leq \beta \leq 1$, fix $z \in \mathbb{R}^d$, we define the Bregman proximal gradient mapping as:
\begin{align}
T_{\lambda}(y):= \arg\min \{ \lambda g(u) + \langle \lambda\nabla f(y), u-y \rangle + \langle u, \beta(z-y) \rangle+ D_h(u, y): u \in \mathbb{R}^d \}.
\end{align}
We could verify that $T_{\lambda}(y)$ is well posed under our settings by the following proposition.
\begin{proposition}\label{Proposition:1}
Suppose that Assumptions \ref{Assumption:A} and \ref{Assumption:B} hold, let $y \in \text{int}~\text{dom}~h$ and $0< \lambda \leq 1/L, 0\leq \beta \leq 1$. Then, the set $T_{\lambda}(y)$ is a nonempty and compact set. 
\end{proposition}

We are now ready to present our Inertial Bregman Proximal Gradient algorithm.
In the algorithm, we require that any step-size $0< \lambda_k \leq 1/L$.  
\begin{algorithm}[htb]
\caption{Inertial Bregman Proximal Gradient - iBPG \label{Alg:1}}
\begin{algorithmic}
\STATE{\textbf{Data:} A function $h$ defined in Definition 1 such that $L$-smad holds on $S$. }
\STATE{\textbf{Initialization:} $x^0=x^{-1} \in \text{int}\text{dom}~h$ and $0< \lambda_k \leq 1/L, 0 \leq \beta_k<1$.}
\STATE{ For $k =0, 1, 2, \dots$, \textbf{repeat}}
\STATE{ \begin{align}
x^{k+1} \in  \arg\min \left\{x: \lambda_k g(x)+ \left\langle x, \lambda_k \nabla f(x^k) - \beta_k(x^k-x^{k-1}) \right\rangle + D_h(x, x^k), x \in \mathbb{R}^d \right\}. \label{Eq:Alg2_2}
\end{align}}
\STATE{ \textbf{until} EXIT received}
\end{algorithmic}
\end{algorithm}

By Proposition \ref{Proposition:1}, under  Assumptions \ref{Assumption:A} and \ref{Assumption:B}, the algorithm is well-defined.
When $h(x)= \frac{1}{2}\|x\|^2$ and consequently $D_h(x, y)=\frac{1}{2}\|x-y\|^2$, the iBPG would be reduced to ipiano algorithm\cite{ochs2014ipiano} for the case where $g(x)$ is convex. When $\beta=0$, iBPG is reduced to BPG algorithm proposed in \cite{bolte2017first}. In this paper, iBPG is a more general inertial nonconvex and nonsmooth algorithm than these two algorithms,  in other words, it can also be treated as an extension of the above two algorithms. 

\section{Convergence Analysis of iBPG}
Throughout the whole analysis of iBPG, we take the following as our blanket assumption
\begin{enumerate}
\item Assumption A and B hold.
\item $h$ is $\sigma$-strongly convex on $\mathbb{R}^d$.
\end{enumerate}

Before analyzing convergence results of iteration seqeunce generated by iBPG for solving problem (P), it is common to show a  general descent(not strictly) lemma firstly in many papers, refer to \cite[Lemma 5]{attouch2013convergence}. Then  along with a similar line, one can analyze the convergence clearly and  succinctly. Such a descent lemma always plays an important role in the whole analysis.
Without exception, we first prove the following lemma, although it can not show the monotone property of objective function $\Psi$ directly. But for a special auxiliary sequence, the monotone property will be present in Lemma \ref{Lemma:2}, which is implied by Lemma \ref{Lemma:1}.
\begin{lemma} \label{Lemma:1}
Let $\{x^k\}_{k\in \mathbb{N}}$ be a sequence generated by iBPG, then
\begin{align}\label{Eq:Lem1_1}
 \Psi(x^{k+1}) +( \lambda_k^{-1} -L - \frac{\beta_k}{\sigma} \lambda_k^{-1}) D_h(x^{k+1},x^k) \leq  \Psi(x^k) + \frac{\beta_k}{\sigma} \lambda_k^{-1}D_h(x^k, x^{k-1}), \forall k \in \mathbb{N}.
\end{align}
\end{lemma}

After adding an inertial term, it can be hard to justify monotonicity of the sequence $\Psi(x^k)$,
seen from \cite{ochs2014ipiano}. We change to give a descent lemma for an auxiliary sequence defined by
$$H_{k,M} := \Psi(x^k) + M D_h(x^{k},x^{k-1}),\forall k \in \mathbb{N}.$$
This descent lemma plays an important role in the following analysis.
In order to guarantee the non-increasing property of the sequences  $\{H_{k,M}\}_{k\in \mathbb{N}}$, we should make some restrictions on the parameter selection. 
We first make the following denotations  
$$ \bar{\beta} := \sup_{k\in \mathbb{N}} \{ \beta_k \}, ~~~~\bar{\lambda} := \sup_{k\in \mathbb{N}} \{ \lambda_k \}, ~~~~\underline{\lambda} := \inf_{k\in \mathbb{N}} \{ \lambda_k \}.$$
We default $\underline{\lambda} >0$. 

\begin{lemma}\label{Lemma:2}
Let $\{x^k \}_{k\in \mathbb{N}}$ be a sequence generated by iBPG. Then we have 
\begin{align}
H_{k+1,M} - H_{k,M} \leq [M -( \lambda_k^{-1} -L - \frac{\beta_k}{\sigma} \lambda_k^{-1})] D_h(x^{k+1},x^k) -(M - \frac{\beta_k}{\sigma} \lambda_k^{-1})D_h(x^k, x^{k-1}), \forall  k \in \mathbb{N}.
\label{Eq:Lem2_1}
\end{align}
Moreover, if we choose fixed $M$ such that
$$ \frac{\bar{\beta}}{\underline{\lambda}}\frac{1}{\sigma} \leq M \leq \frac{1}{\bar{\lambda}} -L -\frac{1}{\sigma} \frac{\bar{\beta}}{\underline{\lambda}},$$
the sequence $\{H_{k,M}\}$ is nonincreasing and convergent for fixed $M$.
\end{lemma}

The next corollary is an obvious result based on Lemma \ref{Lemma:2}, who analyzes boundness of the sequences produced by iBPG.  Since  $H_{k,M}$ is nonincreasing according to Lemma \ref{Lemma:2}, combining with  Assumption \ref{Assumption:B}, it is easy to verify that $\{x^k \}_{k\in \mathbb{N}}$ generated by iBPG is bounded under the same  parameter setting in Lemma \ref{Lemma:2}.
The boundness  would act as a tool in the following analysis.
\begin{corollary}\label{Corollary:1}
Given the same parameter setting in Lemma \ref{Lemma:2}, then the sequence $\{x^k \}_{k\in \mathbb{N}}$ generated by iBPG is bounded.
\end{corollary}

\subsection{Stationary Convergence}\label{SS:SC}
If the parameters are chosen moderate, as in Lemma \ref{Theorem:1}, the range require that $\beta_k$ is strictly smaller than $\frac{\sigma}{2}(\frac{\underline{\lambda}}{\bar{\lambda}}-\underline{\lambda}L)$, then we could get sufficient decrease of the auxiliary sequence  $\{H_{k,M}\}_{k\in \mathbb{N}}$. As a consequence, 
we can bound the sum of Bregman distance between two iteration points generated by iBPG. Due to the strong convexity, we could get that $\lim_{k\rightarrow \infty} \|x^k- x^{k-1}\| =0$ for the sequence  $\{x^k \}_{k\in \mathbb{N}}$ in $\mathbb{R}^d$ by iBPG.
Besides, we could verify under the same condition each limit point is a critical point. 
Assume that  $\{x^k \}_{k\in \mathbb{N}}$ was generated from a starting point $x^0(x^{-1}=x^0)$. 
The set of all limit points of $\{x^k \}_{k\in \mathbb{N}}$ is denoted by
$$ \omega(x^0):= \{ \bar{x}: \text{~an increasing sequence of integers~} \{k_i\}_{i \in \mathbb{N}} \text{~such that~} x^{k_i} \rightarrow \bar{x}  \text{~as ~} i \rightarrow \infty \}.$$
According to the boundness of $\{x^k\}_{k\in \mathbb{N}}$ from Corollary \ref{Corollary:1}, we get $\omega(x^0)$ is nonempty. The following lemma shows that for any $x^0 \in \mathbb{R}^d$, $\omega(x^0) \subseteq \text{crit}~\Psi $ holds.

\begin{lemma}\label{Theorem:1}
Suppose $ \bar{\beta} < \frac{\sigma}{2}(\frac{\underline{\lambda}}{\bar{\lambda}}-\underline{\lambda}L)$. Let $\{x^k\}_{k\in \mathbb{N}} $ be a sequence generated from $x^0$ by iBPG, then 
\begin{enumerate}[(i)]
\item $\sum_{k=0} ^{\infty} D_h(x^{k+1},x^k) < \infty $ and $\lim_{k\rightarrow \infty} \|x^k- x^{k-1}\| =0$.
\item Any cluster point of $\{x^k\}_{k\in \mathbb{N}} $ is a critical point of $\Psi$.
\item $\zeta := \lim_{k\rightarrow \infty} \Psi(x^k)$ exists and $\Psi \equiv \zeta$ on $\Omega$.
\end{enumerate}
\end{lemma}

\begin{remark} In this paper, we choose parameter as
 $$ \bar{\beta} < \frac{\sigma}{2}(\frac{\underline{\lambda}}{\bar{\lambda}}-\underline{\lambda}L).$$
When the stepsize is fixed, it is represented as
$$ \bar{\beta} < \frac{\sigma}{2}(1-\lambda L), $$
which complies with parameter selecting rules for the iPinao(where $h(x)=\frac{1}{2}\|x\|^2 \text{~and~} \sigma =1$)\cite[Table 1]{ochs2018local}. 
\end{remark}

Next, we prove a global $\mathcal{O}(\frac{1}{K})$ convergence rate for $\|x^k-x^{k-1}\|^2$ of the algorithm. In fact, the linear convergence rate can also be get if we add more assumptions, like KL property and concrete KL exponent(calculus of the KL exponent can refer to \cite{li2017calculus}), based on similar deduction\cite[Theorem 6.3]{bolte2017first}.

\begin{corollary}\label{Corollary:2}
Suppose $ \frac{\bar{\beta}}{\underline{\lambda}}\frac{1}{\sigma}  < \frac{1}{\bar{\lambda}} -L -\frac{1}{\sigma} \frac{\bar{\beta}}{\underline{\lambda}}$. Let $\{x^k\}_{k\in \mathbb{N}} $ be a sequence generated from $x^0$ by iBPG, then for all $K \geq 1$, $\min_{1 \leq k \leq K}\|x^k-x^{k-1}\|^2$ converges with a sublinear rate as $\mathcal{O}(\frac{1}{K})$.
\end{corollary}

\subsection{Global Convergence}
In this section, we focus on verifying global convergence results. On the basis of subsection \ref{SS:SC}, we aim to prove that  the sequence $\{x^k\}_{k\in \mathbb{N}} $  which is generated by iBPG converges to a critical point of $\Psi$ defined in (\ref{Eq:P}). In order to prove global convergence, we borrow the proof methodology 
 in \cite{attouch2013convergence}. This proof methodology  prove global convergence result  for several types of nonconvex nonsmooth problem. Other similar forms were referred in many works \cite[Section 3.2]{ochs2014ipiano}\cite[Section 4]{pock2016inertial}\cite[Section 4.2]{bolte2017first}. 
 For the reader’s convenience, we firstly describe the proof methodology summarized in \cite[Theorem 3.7]{ochs2014ipiano} with a few modifications and then apply it to prove convergence of iBPG in  Theorem \ref{Theorem:3}.
\begin{theorem} \cite[Theorem 3.7]{ochs2014ipiano} \label{Theorem:2}
Let $F: \mathbb{R}^{2d} \rightarrow (-\infty, \infty]$ be a proper lower semi-continuous function. Assume that  $\{z^k\}_{k \in \mathbb{N}}:=\{(x^k, x^{k-1})\}_{k \in \mathbb{N}}$ be a sequence generated by a general algorithm from $z^0:=(x^0,x^0)$,  for which the following three conditions are satisfied for any $k \in \mathbb{N}$.
\begin{description}
          \item [(H1)] For each $k \in \mathbb{N}$, there exist a positive $a$ such that 
	                 $$F(z^{k+1}) + a \|x^k-x^{k-1}\|^2 \leq F(z^{k}), \forall k \in \mathbb{N}. $$
	  \item [(H2)] For each $k \in \mathbb{N}$, there exist  a positive $b$ such that for some $v^{k+1} \in \partial F(z^{k+1})$ we have
	  $$\|v^{k+1} \| \leq \frac{b}{2}(\|x^{k+1}-x^{k}\|+ \|x^k-x^{k-1}\|), \forall k \in \mathbb{N}.  $$
    \item [(H3)] There exists a subsequence $(z^{k_j})_{j\in \mathbb{N}}$ such that
    $z^{k_j} \rightarrow \tilde{z}$ and $ F(z^{k_j})\rightarrow F(\tilde{z})$.
\end{description}
Moreover,  if $F$ have the Kurdyka-{\L}􏰀ojasiewicz property at the limit point $\tilde{z}=(\tilde{x},\tilde{x} )$ specified in (H3). Then,  the sequence $\{x^k\}_{k \in \mathbb{N}}$ has finite length, i.e., 􏰒$\sum_{k=1}^{\infty}\|x^k-x^{k-1}\| < \infty$, and converges to $\bar{x} =\tilde{x}$ as $k\rightarrow \infty$, where $(\bar{x}, \bar{x})$ is a critical point of $F$.
 \end{theorem}
 
In our paper, what we need is to verify  conditions in Theorem \ref{Theorem:2} are satisfied
for a sequence $(z^k)_{k\in \mathbb{N}}:= (x^k, x^{k-1})_{k\in \mathbb{N}} \in \mathbb{R}^{2d}$ generated by iBPG. But there is some difference with the convergence  analysis in  \cite{ochs2014ipiano}, we here employ a new function $$H_{M}(x, y) = \Psi(x) + M D_h(x, y)$$  as the original $F$ in Theorem \ref{Theorem:2}. 
Besides, in order to guarantee (H2) holds, we need another assumption(in which the first part was also required in \cite[see Assumption D.(ii)]{bolte2017first}):
\begin{assumption}\label{Assumption:2}
$\nabla h,\nabla f$ are Lipschitz continuous on any bounded subset of $\mathbb{R}^d$, and $\|\nabla^2 h\|$ exists and has an upper bound on any bounded subset of $\mathbb{R}^d$.
\end{assumption} 

The next task is to verifying the three conditions  one by one. Then combining with Theorem  \ref{Theorem:2}, we get the last conclusions that under proper parameter selecting the whole sequence generated by iBPG converge to a critical point. 
\begin{theorem}\label{Theorem:3}.
Suppose $ \frac{\bar{\beta}}{\underline{\lambda}}\frac{1}{\sigma}  < \frac{1}{\bar{\lambda}} -L -\frac{1}{\sigma} \frac{\bar{\beta}}{\underline{\lambda}}$. Let $\{x^k\}_{k\in \mathbb{N}} $ be a sequence generated from $x^0$ by iBPG.  
If $H_{M}(x, y)$ satisfies the Kurdyka--{\L}􏰒ojasiewicz property at some limit point $\tilde{z}=(\tilde{x},\tilde{x})\in \mathbb{R}^{2d}$ and Assumption \ref{Assumption:2} holds.
Then
\begin{enumerate}[(i)]
\item The sequence $\{x^k\}_{k\in \mathbb{N}} $ has finite length, i.e. $\sum_{k=1}^{\infty}\|x^k-x^{k-1}\| <\infty$.
\item $x^k \rightarrow \tilde{x}$ as $k \rightarrow \infty$, and $\tilde{x}$ is a critical point of $\Psi$.
\end{enumerate}
\end{theorem}

\section*{Acknowledgments}
We are grateful for the support from the National Science Foundation of China (No.11501569).



\section*{Appendix}
\subsection*{A  ~Proof of Proposition \ref{Proposition:1}}
Fix any $y, z \in \text{int}~\text{dom}~h$ and $0< \lambda \leq 1/L, 0 \leq \beta \leq 1$. For any $u \in \mathbb{R}^d$, we define 
$$\Psi_h(u) = g(u)+ f(y) + \left\langle u-y ,\nabla f(y)\right\rangle + \lambda^{-1}
\langle u, \beta(z-y) \rangle+\lambda^{-1} D_h(u, y) $$
so that $T_{\lambda}(y) = \arg\min_{u\in \mathbb{R}^d} \Psi_h(u),$
It can also be represented as 
\begin{align*}
\Psi_h(u) &= \Psi(u) -f(u)+f(y)+ \left\langle u-y ,\nabla f(y)\right\rangle + \lambda^{-1} D_h(u,y)+\lambda^{-1}\langle u, \beta(z-y) \rangle \\
& \geq  \Psi(u) + LD_h(u,y)-[f(u)-f(y)- \left\langle u-y ,\nabla f(y)\right\rangle]+ \lambda^{-1}
\langle u, \beta(z-y) \rangle \\
& \geq \Psi(u)+\lambda^{-1}\langle u, \beta(z-y) \rangle.
\end{align*}
where the second inequality is according to $\lambda^{-1} \geq L$ and the last inequality is according to $(f,h)$ is $L$-smooth adaptable. 
Since $\Psi$ is level-bounded, i.e. $\lim_{\|u\| \rightarrow \infty} \Psi(u) = \infty$, there is 
$$ \Psi_h(u) \geq  \|u\|(\frac{ \Psi(u)}{\|u\|}-\lambda^{-1} \|\beta(z-y)\|).$$ 
 Passing to the limit $\|u\| \rightarrow \infty$, the supercoercivity of $\Psi$ implies that $\lim_{\|u\| \rightarrow \infty} \Psi_h(u) = \infty$.
Since $\Psi_h$ is also proper and lower semicontinuous, invoking the modern form of Weierstrass’ theorem (see, e.g., \cite[Theorem 1.9, page 11]{rockafellar2015convex}), it follows that the value $\inf_{\mathbb{R}^d} \Psi_h$ is finite, and the set $\arg\min_{u \in \mathbb{R}^d} \Psi_h(u) \equiv T_{\lambda}(y)$ is nonempty and compact.

\subsection*{B ~Proof of Lemma \ref{Lemma:1} and Lemma \ref{Lemma:2}}
\textbf{Proof of Lemma  \ref{Lemma:1}: }\\
According to the representation (\ref{Eq:Alg2_2}) in iBPG, we have 
\begin{align}\label{Eq:Lem1_2}
\lambda_k g(x^{k+1})+ \left\langle x^{k+1}-x^k, \lambda_k \nabla f(x^k) - \beta_k(x^k-x^{k-1}) \right\rangle + D_h(x^{k+1},x^k) \leq \lambda_k g(x^k), \forall k \in \mathbb{N}.
\end{align}
Besides, according to the $L$-smad property of $(f, h)$, there is
\begin{align}\label{Eq:Lem1_3}
f(x^{k+1}) \leq f(x^k) + \left\langle x^{k+1}-x^k, \nabla f(x^k) \right\rangle + LD_h(x^{k+1},x^k),\forall k \in \mathbb{N}.
\end{align}
Combining (\ref{Eq:Lem1_2}) and (\ref{Eq:Lem1_3}), we get
\begin{align}\label{Eq:Lem1_4}
\lambda_k \Psi(x^{k+1})+ (1-\lambda_k L)D_h(x^{k+1},x^k) 
 \leq \left\langle x^{k+1}-x^k, \beta_k(x^k-x^{k-1}) \right\rangle + \lambda_k g(x^k),\forall k \in \mathbb{N}.
\end{align}
Now, for $ k \geq 0$, using the fact that $\langle p, q \rangle \leq  (1/2) \|p\|^2 + (1/2) \|q\|^2$ for any two vectors $p, q \in \mathbb{R}^d$, yields
\begin{align}
\langle x^{k+1}-x^k, \beta_k(x^k-x^{k-1}) \rangle &\leq (1/2)\beta_k \|x^{k+1}-x^k\|^2 + (1/2)\beta_k \|x^k-x^{k-1}\|^2 \\
&\leq \frac{\beta_k}{\sigma} D_h(x^{k+1},x^k) + \frac{\beta_k}{\sigma}D_h(x^k, x^{k-1}),
\end{align}
where the second inequality is according to the strong convexity of $h$. After rearranging the terms, we get the conclusion. 
\vspace{12pt}

\noindent \textbf{Proof of Lemma  \ref{Lemma:2}:} \\
Clearly, with the definition of $H_{k,M}$,  (\ref{Eq:Lem2_1}) is directly from (\ref{Eq:Lem1_1}).

For the fixed $M$ defined as above, there is 
\begin{align}
\frac{\beta_k}{\sigma} \lambda_k^{-1} \leq \frac{\bar{\beta}}{\underline{\lambda}}\frac{1}{\sigma} \leq M \leq \frac{1}{\bar{\lambda}} -L -\frac{1}{\sigma} \frac{\bar{\beta}}{\underline{\lambda}}  \leq \lambda_k^{-1} -L - \frac{\beta_k}{\sigma} \lambda_k^{-1}, \forall k \in \mathbb{N}.
\end{align}
Consequently, we find that 
$$ H_{k+1,M} \leq H_{k,M},\forall k \in \mathbb{N}. $$ 
Recall that $H_{k,M} \geq \inf \Psi > -\infty$ and $H_{k,M}$ is nonincreasing. This implies the conclusion that $\{H_{k,M}\}$ is convergent .

\subsection*{C  ~Proof of Proof of Lemma \ref{Theorem:1}}

\noindent(i) Since $M \leq \frac{1}{\bar{\lambda}} -L -\frac{s}{\sigma} \frac{\bar{\beta}}{\underline{\lambda}}$, according to  (\ref{Eq:Lem2_1}) and nonnegativeness of $D_h(\cdot,\cdot)$, we have 
\begin{align}
\left(M - \frac{\beta_k}{s\sigma} \lambda_k^{-1}\right)D_h(x^k, x^{k-1} ) \leq  H_{k,M}-H_{k+1,M},\forall k \in \mathbb{N},
\label{Eq:Lem3_1} 
\end{align}
which implies 
\begin{align}
0 \leq \sum_{i=0}^K  \left(M - \frac{\beta_k}{s\sigma} \lambda_k^{-1}\right)D_h(x^k, x^{k-1}) \leq H_{0,M} - H_{K+1,M},
\end{align}
by summing both sides of (\ref{Eq:Lem3_1}) from $0$ to $K$. Since $\{H_{k,M}\}$ is convergent by Lemma \ref{Lemma:2}, letting $K \rightarrow \infty$, we conclude that the infinite sum exists and is finite, i.e.,
$$\sum_{i=0}^K  \left(M - \frac{\beta_k}{s\sigma} \lambda_k^{-1}\right)D_h(x^k, x^{k-1}) < \infty. $$ 
With the definition of $\bar{\beta}$ and $\underline{\lambda}$, we have
\begin{align}
\sum_{i=0}^K \left(M - \frac{\bar{\beta}}{\underline{\lambda}}\frac{1}{s\sigma}\right)D_h(x^k, x^{k-1}) \leq \sum_{i=0}^K  \left(M - \frac{\beta_k}{s\sigma} \lambda_k^{-1}\right)D_h(x^k, x^{k-1}) < \infty. 
\end{align}
If we choose $\bar{\beta}$ such that $ \bar{\beta} < \frac{\sigma}{2}(\frac{\underline{\lambda}}{\bar{\lambda}}-\underline{\lambda}L)$, one can fix $M$ as $ \frac{\bar{\beta}}{\underline{\lambda}}\frac{1}{\sigma} < M \leq \frac{1}{\bar{\lambda}} -L -\frac{1}{\sigma} \frac{\bar{\beta}}{\underline{\lambda}}$. Consequently we get
$ \sum_{k=0} ^{\infty} D_h(x^k, x^{k-1}) < \infty$.
With strong convexity of $h$, there is
$\sum_{k=0} ^{\infty} \|x^k- x^{k-1}\|^2 \leq \frac{2}{\sigma}\sum_{k=0} ^{\infty} D_h(x^k, x^{k-1}) < \infty$ and  thus $\lim_{k\rightarrow \infty} \|x^k- x^{k-1}\| =0$ holds.
\vspace{12pt}

\noindent(ii) Let $\bar{x}$ be an critical point of $\{x^k\}_{k\in \mathbb{N}} $. Let $\{x^{k_i}\}$  be a subsequence such that $\lim_{i\rightarrow \infty} x^{k_i} = \bar{x}$.

By using the first-order optimality condition of the minimization problem (\ref{Eq:Alg2_2}), we obtain
\begin{align*}
0 \in \lambda_{k_i-1} \partial g(x^{k_i}) + \lambda_{k_i-1} \nabla f(x^{k_i-1})- \beta_{k_i-1}(x^{k_i-1}-x^{k_i-2}) + (\nabla h(x^{k_i}) - \nabla h(x^{k_i-1})), \forall k_i \in \mathbb{N}.
\end{align*}
Therefore we obviously observe that
\begin{align}
 \nabla f(x^{k_i}) -\nabla f(x^{k_i-1})  +  \frac{\beta_{k_i-1}}{\lambda_{k_i-1}}(x^{k_i-1}-x^{k_i-2})-\lambda_{k_i-1}^{-1}(\nabla h(x^{k_i}) - \nabla h(x^{k_i-1})) \in \partial \Psi(x^{k_i}), \forall k_i \in \mathbb{N}.  \label{Theorem1_Eq1}
\end{align}
On one hand, limit of the left hand in (\ref{Theorem1_Eq1}) can be justified by the following
\begin{align}
&\|\nabla f(x^{k_i}) -\nabla f(x^{k_i-1})  +  \frac{\beta_{k_i-1}}{\lambda_{k_i-1}}(x^{k_i-1}-x^{k_i-2})-\lambda_{k_i-1}^{-1}(\nabla h(x^{k_i}) - \nabla h(x^{k_i-1}))\| \nonumber\\
&\leq \|\nabla f(x^{k_i}) -\nabla f(x^{k_i-1})\|+ \frac{\bar{\beta}}{\underline{\lambda}}\|x^{k_i-1}-x^{k_i-2}\|+\underline{\lambda}^{-1}\|\nabla h(x^{k_i}) - \nabla h(x^{k_i-1})\| \rightarrow 0, k_i \rightarrow \infty, \label{Theorem1_Eq2}
\end{align}
where the limit can be get according to (i) and the continuity of $\nabla f, \nabla h$. Thus we get 
there exist $v^{k_i} \in  \partial \Psi(x^{k_i})$ such that $\| v^{k_i}\| \rightarrow 0, k_i \rightarrow \infty$.

On the other hand, in view of Lemma \ref{Lemma:2} and (i), the sequence $\{H_{k,M} \}$ is convergent and $D_h(x^k, x^{k-1}) \rightarrow 0$, these together with the definition of $H_{k,M}$ imply that $\lim_{k\rightarrow \infty} \Psi(x^k)$ exists.
Then we derive that $\Psi(x^{k_i}) \rightarrow \Psi(\bar{x}), k_i \rightarrow \infty$. From the lower semicontinuity of $\Psi$, we have 
\begin{align}
 \Psi(\bar{x})  \leq \lim \inf_{i \rightarrow \infty} \Psi(x^{k_i}). \label{Theorem1_Eq3}
\end{align}
According to the iteration step (\ref{Eq:Alg2_2}), for $k_i \geq 1$, we have
\begin{align*}
\lambda_{k_i-1} g(x^{k_i})+ \left\langle x^{k_i}-\bar{x}, \lambda_{k_i-1} \nabla f(x^{k_i-1}) - \beta_{k_i-1}(x^{k_i-1}-x^{k_i-2}) \right\rangle + D_h(x^{k_i},x^{k_i-1}) \leq \lambda_{k_i-1} g(\bar{x})+D_h(\bar{x},x^{k_i-1}) .
\end{align*}
Adding $\lambda_{k_i-1}f(x^{k_i})$ to both sides,  we have 
\begin{align}\label{Theorem1_Eq4}
&\lambda_{k_i-1}\Psi(x^{k_i})+ \left\langle x^{k_i}-\bar{x}, \lambda_{k_i-1} \nabla f(x^{k_i-1}) - \beta_{k_i-1}(x^{k_i-1}-x^{k_i-2}) \right\rangle + D_h(x^{k_i},x^{k_i-1}) \nonumber\\
& \leq \lambda_{k_i-1} g(\bar{x})+\lambda_{k_i-1}f(x^{k_i})+D_h(\bar{x},x^{k_i-1}), \forall k_i \in \mathbb{N} .
\end{align}
After rearranging terms, for all $ k_i \in \mathbb{N}$, there is
\begin{align}\label{Theorem1_Eq5}
\Psi(x^{k_i}) \leq & \Psi(\bar{x})+ f(x^{k_i})-f(\bar{x}) - \left\langle x^{k_i}-\bar{x}, \nabla f(x^{k_i-1}) \right\rangle  + \frac{\beta_{k_i-1}}{\lambda_{k_i-1}}\left\langle x^{k_i}-\bar{x}, x^{k_i-1}-x^{k_i-2} \right\rangle \nonumber\\
& -\lambda_{k_i-1}^{-1}D_h(x^{k_i},x^{k_i-1}) +\lambda_{k_i-1}^{-1}D_h(\bar{x},x^{k_i-1}), \forall k_i \in \mathbb{N}.
\end{align}
 $L$-smad property of $(f,h)$ implies that  for all $ k_i \in \mathbb{N}$
 \begin{align}\label{Theorem1_Eq6}
f(x^{k_i})-f(\bar{x}) - \left\langle x^{k_i}-\bar{x}, \nabla f(x^{k_i-1}) \right\rangle &\leq LD_h(x^{k_i},\bar{x}) +  \left\langle x^{k_i}-\bar{x}, \nabla f(\bar{x})-\nabla f(x^{k_i-1}) \right\rangle \nonumber\\
&= LD_h(x^{k_i},\bar{x}) +  D_f(x^{k_i},x^{k_i-1}) - D_f(x^{k_i},\bar{x})-D_f(\bar{x},x^{k_i-1}). \nonumber\\
&\leq LD_h(x^{k_i},\bar{x}) + L D_h(x^{k_i},x^{k_i-1}) + L D_h(x^{k_i},\bar{x})+L D_h(\bar{x},x^{k_i-1}) 
\end{align}
Besides, (i) implies that
\begin{align}\label{Theorem1_Eq7}
 \frac{\beta_{k_i-1}}{\lambda_{k_i-1}}\left\langle x^{k_i}-\bar{x}, x^{k_i-1}-x^{k_i-2} \right\rangle \leq  \frac{\bar{\beta}}{\underline{\lambda}}\| x^{k_i}-\bar{x}\|\|x^{k_i-1}-x^{k_i-2}\|\rightarrow 0, \forall k_i \in \mathbb{N}.
\end{align}
Plugging (\ref{Theorem1_Eq5}) and (\ref{Theorem1_Eq6}) in (\ref{Theorem1_Eq7}),  passing to the limit, combining with the relationship as $ \underline{\lambda} \leq \lambda_{k_i} \leq \bar{\lambda}$,  we have
\begin{align}\label{Theorem1_Eq8}
\lim_{i\rightarrow \infty} \Psi(x^{k_i}) &\leq  \Psi(\bar{x})+ \lim_{i\rightarrow \infty}\left[(-\bar{\lambda}^{-1}+L)D_h(x^{k_i},x^{k_i-1}) + (\underline{\lambda}^{-1}+L)D_h(\bar{x},x^{k_i-1})+2L D_h(x^{k_i}, \bar{x})\right] \nonumber\\
 &\leq  \Psi(\bar{x})+ \lim_{i\rightarrow \infty} (\underline{\lambda}^{-1}+L)\left[D_h(\bar{x},x^{k_i-1})+D_h(x^{k_i}, \bar{x})\right],
\end{align}
where the second inequality is based on $L \leq \bar{\lambda}^{-1} \leq \underline{\lambda}^{-1}$ in iBPG. From Lemma \ref{Lemma:1}(i), $ D_h(x^{k_i},x^{k_i-1})\rightarrow 0$. 
Combining with the continuity of $\nabla h$, we have
\begin{align*}
 \lim_{i\rightarrow \infty} \left[D_h(\bar{x},x^{k_i-1})+D_h(x^{k_i}, \bar{x})\right]
&= \lim_{i\rightarrow \infty} \left[D_h(x^{k_i},x^{k_i-1}) + \langle \nabla h(x^{k_i-1})-\nabla h(\bar{x}) , x^{k_i}-\bar{x}\rangle \right] \nonumber\\
&\leq  \lim_{i\rightarrow \infty} \left[D_h(x^{k_i},x^{k_i-1}) + \|\nabla h(x^{k_i-1})-\nabla h(\bar{x})\|\| x^{k_i}-\bar{x}\| \right] \nonumber\\
&= 0.
\end{align*}
Hence we have
\begin{align}\label{Theorem1_Eq9}
\lim \sup_{i\rightarrow \infty} \Psi(x^{k_i}) \leq \Psi(\bar{x}).
\end{align}
Combining (\ref{Theorem1_Eq3}) and (\ref{Theorem1_Eq9}) yields  $\Psi(x^{k_i}) \rightarrow \Psi(\bar{x}), k_i \rightarrow \infty$. 

Thus, according to these two hand,  and  the closedness of $\partial \Psi$ (see, Definition 1), we have $0 \in \partial \Psi(\bar{x})$.
\vspace{12pt}

\noindent(iii)  According to the last part of the proof in (ii), we know that if $\omega(x^0) =\emptyset$, the conclusion holds;
Otherwise, take $\bar{x} \in \omega(x^0)$ with a convergent subsequences $\{x^{k_i}\}$ meet that $\lim_{i\rightarrow \infty} x^{k_i} = \bar{x}$, there is
 $$\zeta = \lim_{i\rightarrow \infty} \Psi(x^{k_i}) = \Psi(\bar{x}).$$
Thus the conclusion is completed since $\bar{x}$ is arbitrary.

\subsection*{D ~Proof of Corollary \ref{Corollary:2}}

Setting $M = \frac{1}{\bar{\lambda}} -L -\frac{1}{\sigma} \frac{\bar{\beta}}{\underline{\lambda}}$, we get the following from Lemma \ref{Lemma:2}
\begin{align}
 \frac{\sigma}{2}\left(M - \frac{\bar{\beta}}{\underline{\lambda}}\frac{1}{\sigma} \right)\sum_{i=1}^K \|x^k-x^{k-1}\|^2   
&\leq \sum_{i=1}^K \left(M - \frac{\bar{\beta}}{\underline{\lambda}}\frac{1}{\sigma} \right)D_h(x^k, x^{k-1})  \nonumber \\
&\leq \sum_{i=1}^K  \left(M - \frac{\beta_k}{\sigma} \lambda_k^{-1}\right)D_h(x^k, x^{k-1}) \nonumber \\
&\leq H_{1,M} - H_{K+1,M}.\nonumber
\end{align}
Thus we get
$$\min_{1 \leq k \leq K}\|x^k-x^{k-1}\|^2 \leq \frac{H_{1,M}-H_{K+1,M}}{K (\frac{\sigma}{2}M-\frac{\bar{\beta}}{\underline{\lambda}})} \leq \frac{H_{1,M}-\Psi^{*}}{K (\frac{\sigma}{2}M-\frac{\bar{\beta}}{\underline{\lambda}})},$$  
the proof is concluded.

\subsection*{E ~Proof of Theorem \ref{Theorem:3}} 

According to Theorem \ref{Theorem:2}, combining three conditions illustrated in Theorem  \ref{Theorem:2} and KL property at $\tilde{z}$ could guarantee that conclusion (i) holds. Conclusion (ii) are followed by Theorem \ref{Theorem:3}(i).

Next we will verify the three conditions for iBPG.

(i) From Lemma \ref{Lemma:2}, since $h$ is strong convex,  one can show inductively that (H1) holds.
\vspace{6pt}

(ii)From (\ref{Eq:Alg2_2}), there exists $v_{k+1} \in \partial H_{M}(x^{k+1},x^k)$ such that
\begin{align*}
v_{k+1} = \left(\nabla f(x^{k+1})-\nabla f(x^{k}) + \frac{\beta_{k+1}}{\lambda_{k+1}}(x^k-x^{k-1})+ (M-\lambda_{k+1}^{-1}) (\nabla h(x^{k+1})- \nabla h(x^k)), \langle -\nabla^2 h(x^{k}), x^{k+1} -x^k \rangle \right).
\end{align*}
Due to Corollary \ref{Corollary:1}, $\{x^k\}_{k\in \mathbb{N}}$ generated by iBPG is bounded. Thus according to Assumption \ref{Assumption:2}, there exist 
 $L_f,L_h$ such that for any $k \in \mathbb{N}$, $\|\nabla h(x^{k+1}) - \nabla h(x^{k})\| \leq L_h\|x^{k+1}-x^{k}\|, \|\nabla f(x^{k+1}) - \nabla f(x^{k})\| \leq L_f\|x^{k+1}-x^{k}\|$.
Combining $\|\nabla^2 h(x)\|$ has an upper bound(set as $\delta$)  on $\{x^k\}_{k\in \mathbb{N}}$, we have
\begin{align}
 \|v_{k+1}\| &\leq \left(L_f+ (M-\lambda_{k+1}^{-1})L_h+\|\nabla^2 h(x^k)\| \right)\|x^{k+1}-x^{k}\| +  \frac{\beta_{k+1}}{\lambda_{k+1}} \|x^k-x^{k-1}\| \nonumber \\
&\leq \left(L_f+ (M-\bar{\lambda}^{-1})L_h+ \delta \right)\|x^{k+1}-x^{k}\| +  \frac{\bar{\beta}}{\underline{\lambda}} \|x^k-x^{k-1}\|
\end{align}
Denote $b:= \max\{L_f+ (M-\bar{\lambda}^{-1})L_h+ \delta ,  \frac{\bar{\beta}}{\underline{\lambda}} \}$, (H2) is satisfied.
\vspace{6pt}

(iii)(H3) naturally follows from Lemma \ref{Theorem:1}.

\end{document}